# Some remarks about Hankel determinants which are related to Catalan-like numbers

Johann Cigler

**Abstract.**

This note presents some results about Hankel determinants of backwards shifted Catalan-like numbers and related sequences.

**1. Introduction.**

Let $\mathbf{s} = (s_k)_{k \geq 0}$ be a sequence of real or complex numbers. Define numbers $a_{n,k} = a_{n,k}(\mathbf{s})$ by

$$a_{n,k} = a_{n-1,k-1} + s_k a_{n-1,k} + a_{n-1,k+1} \tag{1}$$

with $a_{n,k} = 0$ for $k < 0$ and $a_{0,k} = [k = 0]$. We also extend $a_{n,k}$ to $n \in \mathbb{Z}$ by setting $a_{n,k} = 0$ for $n < 0$.

Following Aigner [1] we call the matrices $(a_{n,k}(\mathbf{s}))$ *admissible matrices* and the numbers $a_n(\mathbf{s}) = a_{n,0}(\mathbf{s})$ *Catalan-like numbers of type* $\mathbf{s}$.

Let

$$D_{m,k,n}(\mathbf{s}) = \det\left(a_{i+j+m,k}(\mathbf{s})\right)_{i,j=0}^{n-1} \tag{2}$$

with $D_{m,k,0}(\mathbf{s}) = 1$. For $k = 0$ we also write briefly $D_{m,k,n}(\mathbf{s}) = D_{m,n}(\mathbf{s})$.

It is well known (cf. [1]) that

$$\begin{aligned} D_{0,n}(\mathbf{s}) &= 1, \\ D_{1,n}(\mathbf{s}) &= s_{n-1} D_{1,n-1}(\mathbf{s}) - D_{1,n-2}(\mathbf{s}). \end{aligned} \tag{3}$$

with initial values $D_{1,0}(\mathbf{s}) = 1$ and $D_{1,1}(\mathbf{s}) = s_0$.

We shall prove the following results for Hankel determinants $D_{m,k,n}(\mathbf{s})$ with negative $m$, which have been conjectured in [4]-[6].

**Theorem 1**

*For a positive integer $m$ we get $D_{-m,n}(\mathbf{s}) = 0$ for $0 < n \leq m$ and*

$$D_{-m,n+m+1}(\mathbf{s}) = (-1)^{\binom{m+1}{2}} D_{m,n}(E\mathbf{s}) \tag{4}$$

*for $n \in \mathbb{N}$, where $E\mathbf{s} = (s_{k+1})_{k \geq 0}$.*

**Remark**

Since $D_{0,n}(\mathbf{s}) = 1$ identity (4) remains true for $m = 0$.

For example, for $\mathbf{s} = (1,0,0,0,\cdots)$ we have $a_n(\mathbf{s}) = \binom{n}{\lfloor \frac{n}{2} \rfloor}$ and

$$(D_{2,0,n}(\mathbf{s}))_{n \geq 0} = (1,2,3,4,5,6,7,8,9,10,11,12,\cdots),$$

$$(D_{2,0,n}(E\mathbf{s}))_{n \geq 0} = (1,1,2,2,3,3,4,4,5,5,6,6,\cdots), \text{ and}$$

$$(D_{-2,0,n}(\mathbf{s}))_{n \geq 0} = (1,0,0,-1,-1,-2,-2,-3,-3,-4,-4,-5,-5,-6,-6,-7,-7,-8,-8,\cdots).$$

For constant sequences $\mathbf{s} = \underline{c} := (c,c,c,\cdots)$ we also obtain results for $k > 0$.

**Theorem 2**

For $k, m \in \mathbb{N}$ we get $D_{-m,k,n}(\underline{c}) = 0$ for $0 < n \leq m + k$ and

$$D_{-m,k,n+m+k+1}(\underline{c}) = (-1)^{\binom{m+k+1}{2}} D_{m,k,n}(\underline{c}) \tag{5}$$

for $n \in \mathbb{N}$.

For example, for $c = 1$ we get the Motzkin numbers $a_n(\underline{1}) = a_{n,0}(\underline{1}) = M_n$ (cf. OEIS [12], A001006). (5) gives for example,

$$(D_{1,1,n}(\underline{1}))_{n \geq 0} = (1,1,1,1,0,0,-1,-1,-1,-1,0,0,\cdots),$$
$$(D_{-1,1,n}(\underline{1}))_{n \geq 0} = (1,0,0,-1,-1,-1,-1,0,0,1,1,1,1,\cdots),$$
$$(D_{2,2,n}(\underline{1}))_{n \geq 0} = (1,1,0,-4,-4,0,9,9,0,-16,-16,0,25,25,\cdots),$$
$$(D_{-2,2,n}(\underline{1}))_{n \geq 0} = (1,0,0,0,0,1,1,0,-4,-4,0,9,9,\cdots).$$

For $m = 0$ Theorem 2 gives

$$D_{0,k,n+k+1}(\underline{c}) = (-1)^{\binom{k+1}{2}} D_{0,k,n}(\underline{c}) \text{ with } D_{0,k,0}(\underline{c}) = 1 \text{ and } D_{0,k,n}(\underline{c}) = 0 \text{ for } 0 < n \leq k.$$

This implies

**Corollary**

For each $k \in \mathbb{N}$ and each $c$

$$D_{0,k,(k+1)n}(\underline{c}) = (-1)^{\binom{k+1}{2}n}, \tag{6}$$
$$D_{0,k,n}(\underline{c}) = 0 \quad \text{else.}$$

This result is a special case of [8], Theorem 1 by choosing there $x, y$ with $x + y = c$ and $xy = 1$.

**Remark**

It is well known (cf. [11], (5.40), (5.41)) that

$$D_{0,0,n}(c) = 1,$$
$$D_{1,0,n}(c) = F_{n+1}(c), \quad (7)$$
$$D_{2,0,n}(c) = \sum_{j=0}^{n} F_{j+1}(c)^2$$

if $F_n(x)$ denotes the Fibonacci polynomials $F_n(x) = xF_{n-1}(x) - F_{n-2}(x)$ with $F_0(x) = 0$ and $F_1(x) = 1$. However, in [11] it is formulated in different terminology. To see the connection observe that (14) gives for the generating function $A(x,c) = \sum_{n \geq 0} a_n(\underline{c}) x^n$ the continued fraction

$$A(x,c) = \cfrac{1}{1 - cx - \cfrac{x^2}{1 - cx - \cfrac{x^2}{1 - cx - \cdots}}}.$$

Identity (6) is an analog for $k > 0$ of $D_{0,0,n}(c) = 1$. For the sake of completeness, let us mention analogs of the other parts of (7).

Let $L_n(x)$ denote the Lucas polynomials $L_n(x) = xL_{n-1}(x) - L_{n-2}(x)$ with $L_0(x) = 2$ and $L_1(x) = x$. Then [8], Theorem 2, gives

$$D_{1,k,(k+1)n}(c)(-1)^{n\binom{k+1}{2}} = D_{1,k,(k+1)n+k}(c)(-1)^{n\binom{k+1}{2}+\binom{k}{2}} = \frac{y^{(k+1)(n+1)} - x^{(k+1)(n+1)}}{y^{k+1} - x^{k+1}} \quad \text{for } x + y = c \text{ and}$$

$xy = 1$. This implies

$$D_{1,k,(k+1)n}(c)(-1)^{n\binom{k+1}{2}} = D_{1,k,(k+1)n+k}(c)(-1)^{n\binom{k+1}{2}+\binom{k}{2}} = F_{n+1}(L_{k+1}(c)), \quad (8)$$
$$D_{1,k,n}(c) = 0 \text{ else.}$$

To see this, observe that Binet's formulae give $F_n(c) = \dfrac{y^n - x^n}{y - x}$ and $L_k(c) = x^k + y^k$.

For $m = 2$ we guess that for $k \geq 1$

$$D_{2,k,(k+1)n}(c) = (-1)^{n\binom{k+1}{2}} F_{n+1}\left(L_{k+1}(c)\right)^2,$$

$$D_{2,k,(k+1)n+k-1}(c) = (-1)^{\binom{k-1}{2}} D_{2,k,(k+1)n}(c),$$

$$D_{2,k,(k+1)n+k}(c)(-1)^{n\binom{k+1}{2}+\binom{k}{2}} = (k+1)F_{k+1}(c)\sum_{j=0}^{n} F_{j+1}\left(L_{k+1}(c)\right)^2,$$

$$D_{2,k,n}(c) = 0 \quad \text{else.}$$

(9)

It also seems that for $k \geq m-1$

$$D_{m,k,(k+1)n}(c) = (-1)^{n\binom{k+1}{2}} F_{n+1}\left(L_{k+1}(c)\right)^m. \tag{10}$$

## 2. Proofs

Let us recall a well-known combinatorial interpretation of $a_{n,k}(\mathbf{s})$ as the weight of certain lattice paths.

A lattice path which never goes below the $x-$axis, which starts at $(0,0)$ and consists of up-steps $U = (1,1)$, down-steps $D = (1,-1)$ and horizontal steps $H = (1,0)$ will be called a 3-step path. We define the weight $w(P) = w(\mathbf{s}, P)$ of a 3-step path $P$ as the product of the weights of its steps where $w(U) = w(D) = 1$, $w(H_{i,j}) = w((i,j) \to (i+1,j)) = s_j$ and define the weight of a set of paths as the sum of their weights.

Then (1) shows that $a_{n,k}(\mathbf{s})$ is the weight of the set of all 3-step paths from $(0,0)$ to $(n,k)$. Since $a_{0,0}(\mathbf{s}) = 1$ we introduce a trivial path of length 0 with weight 1. A 3-step path which ends on the $x-$axis is called Motzkin path.

The generating function $A(x,\mathbf{s}) = \sum_{n\geq 0} a_n(\mathbf{s})x^n$ of the weight of all Motzkin paths satisfies

$$A(x,\mathbf{s}) = 1 + s_0 x A(x,\mathbf{s}) + x^2 A(x, E\mathbf{s}) A(x,\mathbf{s}). \tag{11}$$

For the trivial path has weight 1, the paths starting with a horizontal step give the second term and paths which start with an up-step give the last term since each such path $P$ can be uniquely written as $P = UP_1 DP_0$, where $P_1$ is a Motzkin path on height 1 and $P_0$ a Motzkin path on height 0.

Identity (11) can also be stated as

$$A(x,\mathbf{s}) = \sum_{n\geq 0} a_n(\mathbf{s})x^n = \frac{1}{1 - s_0 x - x^2 A(x, E\mathbf{s})}. \tag{12}$$

Our proofs depend on the following

**Lemma** ([2] Theorem 2, [3] Proposition 2.5)

*Let* $u(x) = \sum_{n \geq 0} u_n x^n$ *with* $u_0 = 1$ *and* $v(x) = \dfrac{1}{u(x)} = \sum_{n \geq 0} v_n x^n$.

*Setting* $u_n = v_n = 0$ *for* $n < 0$ *we get for* $M \in \mathbb{N}$

$$\det\left(u_{i+j-M}\right)_{i,j=0}^{N+M} = (-1)^{N+\binom{M+1}{2}} \det\left(v_{i+j+M+2}\right)_{i,j=0}^{N-1}. \tag{13}$$

**Proof of Theorem 1**

Choosing $u(x) = A(x, \mathbf{s})$ we get $v(x) = 1 - s_0 x - x^2 A(x, E\mathbf{s})$ by (12). Thus, we get $u_{i+j-M} = a_{i+j-M}(\mathbf{s})$ and $v_{i+j+M+2} = -a_{i+j+M}(E\mathbf{s})$. This implies for $m \in \mathbb{N}$

$$\det\left(a_{i+j-m}(\mathbf{s})\right)_{i,j=0}^{n+m} = (-1)^{n+\binom{m+1}{2}} \det\left(-a_{i+j+m}(E\mathbf{s})\right)_{i,j=0}^{n-1} = (-1)^{\binom{m+1}{2}} \det\left(a_{i+j+m}(E\mathbf{s})\right)_{i,j=0}^{n-1}.$$

Another proof has been given in [6].

**Proof of Theorem 2.**

Let $A(x, \mathbf{c}) = \sum_{n \geq 0} a_{n,0}(\underline{c}) x^n = \sum_{n \geq 0} a_n(\underline{c}) x^n$.

From (11) we get

$$A(x, c) = 1 + cxA(x, c) + x^2 A(x, c)^2 \tag{14}$$

and therefore

$$x^k A(x, c)^{k+1} = x\left(x^{k-1} A(x, c)^k + c\left(x^k A(x, c)^{k+1}\right) + \left(x^{k+1} A(x, c)^{k+2}\right)\right). \tag{15}$$

Comparing coefficients with (1) gives

$$x^k A(x, c)^{k+1} = \sum_{n \geq 0} a_{n,k}(\underline{c}) x^n. \tag{16}$$

Formula (14) implies the closed formula

$$A(x, c) = \frac{1 - cx - \sqrt{(1-cx)^2 - 4x^2}}{2x^2}. \tag{17}$$

These results are well known (cf. [1]).

For the computation of $\dfrac{1}{A(x,c)^k}$ we recall that the bivariate Lucas polynomials

$$L_n(x, s) = xL_{n-1}(x, s) + sL_{n-2}(x, s) \tag{18}$$

with initial values $L_0(x, s) = 2$ and $L_1(x, s) = x$ satisfy

$$L_n(x + y, -xy) = x^n + y^n. \tag{19}$$

To show this it suffices to verify (19) for $n=0$ and $n=1$ and observe that
$$x^n + y^n = (x+y)(x^{n-1} + y^{n-1}) - xy(x^{n-2} + y^{n-2}).$$

By (12) we get $\dfrac{1}{A(x,c)} + x^2 A(x,c) = 1 - cx = L_1(1-cx, -x^2)$ and thus by (19)

$$\frac{1}{A(x,c)^{k+1}} + x^{2k+2} A(x,c)^{k+1} = L_{k+1}(1-cx, -x^2). \tag{20}$$

It should be noted that (20) is closely related to Binet's formula $L_n(x,s) = \alpha^n + \beta^n$ with

$\alpha(x,s) = \dfrac{x - \sqrt{x^2 + 4s}}{2}$ and $\beta(x,s) = \dfrac{x + \sqrt{x^2 + 4s}}{2}$ by observing that

$x^2 A(x,c) = \alpha(1-cx, -x^2)$ and $\dfrac{1}{A(x,c)} = \beta(1-cx, -x^2).$

Choosing $u(x) = A(x,c)^{k+1} = \sum\limits_{n \geq 0} a_{n+k,k}(\underline{c}) x^n$ and $M = m+k$ we get $u_{i+j-M} = a_{i+j-m,k}(\underline{c})$.

Since $L_k(1-cx, -x^2)$ is a polynomial of degree $k$ in $x$ we get

$$v(x) = L_{k+1}(1-cx, -x^2) - x^{2k+2} A(x,c)^{k+1}.$$

This gives $v_{i+j+M+2} = v_{i+j+m+k+2} = -a_{i+j+m,k}(\underline{c})$. For $N = n$ we get (5).

**Remark**

For $c = 0$ we get $A(x, 0) = \dfrac{1 - \sqrt{1-4x^2}}{2x^2} = \sum\limits_{n \geq 0} C_n x^{2n}$ with the Catalan numbers $C_n$. This case has already been considered in [7] with the same method and in [9] with bijective proofs. The Corollary could also be proved by the method of Hankel continued fractions in [10].

Finally let us consider the Hankel determinants of the coefficients of $\dfrac{1}{A(x,c)^{k+1}} = \sum\limits_{n \geq 0} b_{n,k}(\underline{c}) x^n$.

**Theorem 3**

$$d_{0,k,n+1}(\underline{c}) := \det\bigl(b_{i+j,k}(\underline{c})\bigr)_{i,j=0}^{n} = (-1)^n D_{k+2,k,n}(\underline{c}). \tag{21}$$

For $k = 2$ and $c = 1$ we get

$(b_{n,2}(\underline{1}))_{n \geq 0} = (1, -3, 0, 2, 0, 0, -1, -3, -9, -25, -69, -189, -518, \cdots),$

$(a_{n,2}(\underline{1}))_{n \geq 0} = (0, 0, 1, 3, 9, 25, 69, 189, 518, \cdots),$

$(d_{0,2,n}(\underline{1}))_{n \geq 0} = (1, 1, -9, -4, 20, -225, -45, 126, \cdots),$

$(D_{4,2,n}(\underline{1}))_{n \geq 0} = (1, 9, -4, -20, -225, 45, 126, \cdots).$

For example,

$$d_{0,2,3}(\underline{1}) = \det\begin{pmatrix} 1 & -3 & 0 \\ -3 & 0 & 2 \\ 0 & 2 & 0 \end{pmatrix} = -4 = \det\begin{pmatrix} 9 & 25 \\ 25 & 69 \end{pmatrix} = D_{4,2,2}(\underline{1}).$$

**Proof**

Choosing in the Lemma $u(x) = \dfrac{1}{A(x,c)^{k+1}}$, $N = n$, $M = 0$ we get $u_{i+j} = b_{i+j,k}(\underline{c})$ and $v(x) = A(x,c)^{k+1} = \sum_{n \geq 0} a_{n+k,k}(\underline{c})x^n$ and $v_{i+j+M+2} = a_{i+j+k+2,k}(\underline{c})$.

**References**

[1] Martin Aigner, Catalan-like numbers and determinants, J. Comb. Th. A 87 (1999), 33-51

[2] George Andrews and Jet Wimp, Some q-orthogonal polynomials and related Hankel determinants, Rocky Mountain J. Math. 32(2), 2002, 429-442

[3] Johann Cigler, A special class of Hankel determinants, arXiv:1302.4235

[4] Johann Cigler, Shifted Hankel determinants of Catalan numbers and related results II: Backward shifts, arxiv: 2306.07733

[5] Johann Cigler, Some experimental observations about Hankel determinants of convolution powers of Catalan numbers, arXiv: 2308.07642

[6] Johann Cigler, Some results and conjectures about Hankel determinants of sequences which are related to Catalan-like numbers, arXiv:2309.15557

[7] Johann Cigler, Hankel determinants of convolution powers of Catalan numbers revisited, arXiv:2403.11244

[8] Johann Cigler and Christian Krattenthaler, Some determinants of path generating functions, Adv. Appl. Math. 46 (2011), 144-174

[9] Markus Fulmek, Hankel determinants of convoluted Catalan numbers and nonintersecting lattice paths: A bijective proof of Cigler´s conjecture. arXiv:2402.19127

[10] GuoNiu Han, Hankel continued fraction and its applications, arXiv:1406.1593

[11] Christian Krattenthaler, Advanced determinant calculus: a complement, Linear Algebra Appl. 411 (2005), 68-166

[12] OEIS, https://oeis.org/